\documentclass[12pt]{article}
\usepackage{amsfonts}
\usepackage{bbm}
\usepackage{graphicx}
\usepackage{multirow}

\title{Combining Survival Trials Using Aggregate Data Based on Misspecified Models}
\date{\today}
\author{Tinghui Yu \footnote{FDA, Center for Devices and Radiological Health.}, Yabing Mai \footnote{Merck Research Laboratories.}, Sherry Liu $^*$, Xiaofei Hu $^\dagger$}

\begin{document}
\maketitle
\noindent
\begin{center}
\textbf{Abstract}
\end{center}

\noindent
The treatment effects of the same therapy observed from multiple clinical trials can be very different. Yet the patient characteristics accounting for the differences may not be identifiable in real practice so that it is necessary to estimate and report the overall treatment effect for the general popoulation during the development and validation of a new therapy. The non-linear structure of the maximum partial likelihood estimates for the (log) hazard ratio defined with a Cox proportional hazard model leads to challenges in the statistical analyses for combining such clinical trials.  In this paper, we formulated the expected overall treatment effects using various modeling assumptions. Then we proceeded to propose efficient estimates together with a version of Wald test for the combined hazard ratio using only aggregate data. Interpretation of the methods are provided in the framework of robust data analyses involving misspecified models. \\


\noindent
\section{Introduction}
Multiple clinical trials may be performed to validate a newly developed therapy to account for the variability of the targeted patient population. The data collected from each individual clinical trial may be used to generate the efficacy (or safety) estimates for each specific trial as well as the overall efficacy estimate for the overall population to support the effectiveness of the therapy. In this paper, we will focus on the analyses of time-to-event data ($t$) using Cox proportional hazard models. Let $\mathbf{Z}=\{Z_1, \ldots, Z_k\}$ be the covariates available for each enrolled patient. Given $\mathbf{Z}$, the proportional hazard model assumes that the hazard function for a patient from the $i$-th ($i=1,\ldots,M$) trial can be writtens as
\[ h_i(t|\mathbf{Z})=h_{i0}(t)\exp(\tilde{\beta}_i'\mathbf{Z}), \]
where $h_{i0}(t)$ is the baseline hazard function of the $i$-th trial with unknown formulation, $\tilde{\beta}_i=\{\beta_{i1}, \ldots, \beta_{ik}\}$ is the vector of log hazard ratios defined specifically for the patients enrolled in the $i$-th trial.\\

\noindent
Efficient estimates of the trial specific $\tilde{\beta}_i$ based on maximal partial likehood (MPLE) were well developed since the days of Cox (1972, 1975). In this paper, we are more interested in the methods for statistical inference based on the information contained in the pooled data from all $M$ trials. Meta analyses were invented to address such needs. Yuan (2009) provided a comprehensive review of the statistical methods used to generate an overall (log) hazard ratio estimate based on multiple clinical trials. The convex combination of all the log hazard ratio estimates from each individual trial is a popular choice (Wei, Lin and Weissfeld 1989):
\begin{equation}
	\hat{\beta}_{Lj}=\sum_{i=1}^M w_i\hat{\beta}_{ij}, \forall w_1+\ldots+w_M=1.
	\label{linear_estimate}
\end{equation}
Here $\hat{\beta}_{ij}$ ($j=1,\ldots,k$) denotes the MPLEs of ${\beta}_{ij}$ derived from the patient data collected from the $i$-th trial. The weights $w_i$ can be arbitrary constants sum up to unity, among which the inverse variance scheme may be the most popular option due to the obvious advantage of minimized variance among the family of all convex combinations of $\hat{\beta}_{ij}$'s. Though less common, some researchers also recommend using the linear combination of the hazard ratio estimates for any given covariate value $Z_j$:
\begin{equation}
	\hat{b}_{Lj}=\sum_{i=1}^M w_i e^{\hat{\beta}_{i,j}Z_j}, \forall w_1+\ldots+w_M=1.
	\label{linear_estimate_hr}
\end{equation}

 \noindent
To account for the differences in each $\hat{\beta}_{ij}$ due to randomness inherited from the observed data, DerSimonian and Laird (1986) provided an overall treatment effect estimate based on random effect models. Yuan continued to propose a meta-ANOVA model and a meta-polynomial model to address the differences in $\hat{\beta}_{ij}$'s due to mis-match of baseline information from different trials. These methods require a large number of trials ($M$) to support the regression algorithms. Here we are only interested in the cases where very limited number ($M\geq 2$) of trials are deliverable.\\

\noindent
As far as we know from the research literature, all the statistical models developed for meta-analyses assume that there exists a set of ``true'' baseline hazard function $h_{0i}(\cdot)$ as well as a unique ``true'' log hazard ratio (vector) $\tilde{\beta}$ through out all the trials designed for the same therapy. The baseline $h_{0i}(\cdot)$ may vary by $i$, while $\tilde{\beta}$ has to be consistent for all the $M$ trials. The differences in the estimated log hazard ratios $\hat{\beta}_{i}=\{\hat{\beta}_{i1}, \ldots, \hat{\beta}_{ik}\}$ are either attributed to the randomness of the observed outcomes or the incompleteness of the covariate data. However, in some cases, such assumptions may be far from the truth. The treatment effect on different populations may be essentially different due to discrepancies in some latent patient characteristics. The problem is even more acute if the therapy under investigation targets specific gene expressions. The human genome is overwhelmingly complex such that even the most well devised targeted treatment can be subject to unexpected impacts from genes outside the targeted region. Hence the true values of $\tilde{\beta}_i$'s may vary with $i$ because the patient populations enrolled for different trials are in fact heterogeneous with regard to their responses to the therapy. Even though the inclusion/exclusion criteria and the design of these trials may appear to be perfectly aligned, it cannot help to suppress the differences in the patients recruited for different trials. It is impossible to adjust the hazard ratio estimates for these latent controlling factors because they are usually unknown to the researchers or not able to be detected by any currently available technology/assay. Instead, the researcher has to compromise with an overall treatment effect estimate derived from the pooled data, or from the summary statistics reported for each trial if patient line data are not available.\\

\noindent
Mehrotra et. al. (2012) investigated the effect of deviation from the constant hazard ratio assumption underlying the standard stratified survival analyses. Our method proposed in this paper can be used to establish an efficient overall treatment effect estimate covering multiple stratums of high variability. Denne et. al. (2014) and Li (2014) provided another good example illuminating the necessity of considering different treatment effects ($\tilde{\beta}_i$) for different trials. One version of an assay may be used to enroll patients by a biomarker for a validation trial. Another version of the assay targeting the same biomarker may be needed when the treatment is ready for marketing due to advances of technology or solely to reduce the cost and time for patient screening. The two versions of the assay may not perfectly match with each other. The discrepancies in the test results may reflect some differences in the patient's biological profile which are not intended to be captured by either assay. Such differences may affect the patient responses to the therapy. Hence the intent-to-treat patient population is divided into multiple subgroups with different assay result combinations. It is reasonable to assume that the patient responses to the same targeted treatment are in fact different (or even opposite as was observed in some real life examples) across these subgroups defined by both assays (Sargent et. al. 2005).  However, only the market-ready version of the assays will be available to the patients/client laboratories so that it is the overall treatment effect covering all subgroups, rather than the treatment effect on patients with specific test result combinations, that concerns the developer.\\

\noindent
The linear estimate $\hat{\beta}_L$ defined by (\ref{linear_estimate}) or $\hat{b}_L$ by (\ref{linear_estimate_hr}) was used as the estimator for the overall treatment effect in these mentioned publications. The linear estimates has its own merit. It is easy to be implemented for calculation and the derivation of its variance is straightforward. It is guaranteed to be close to the correct answer if the true $\tilde{\beta}_i$ values across different trials only differ by a small amount of no practical significance. However, if the $\tilde{\beta}_i$'s are very different, despite the choice of the weights ($w_i$), the usage of a linear estimate is not mathematically justifiable because each component of the linear combination converges to a completely different $\beta_{ij}$. It is not very likely that the limit of the estimator ($\sum w_i\hat{\beta}_{ij}$) is an exact measurement of the overall treatment effect because the (log) hazard ratio defined by the Cox proportional hazard model is obviously nonlinear with regard to the observed survival times. As a matter of fact, we will show that in most cases $\hat{\beta}_L$ (and $\hat{b}_L$) defined by the inverse variances coefficients or by the proportion of trials sizes leads to over estimate of the true overall hazard with the treatment. \\

\noindent
The major concerns and challenges for a pharmaceutical researcher facing data collected from multiple trials can be summarized by two questions. First, what is the proper definition of the overall treatment effect? Hazard ratio is very different from the other commonly used endpoints such as the mean/median survival time and response rate. The later resorts to a natural measurement (in most cases, counting) of observed events or time, while hazard ratio is an artificial concept specifically invented for the Cox proportional hazard model.  The definition of a ``hazard ratio for the overall population'' is ambiguous now that the Cox model is no more applicable to a mixed population. In the ideal case one may have all the exact knowledge about the true underlying performance of the therapy including but not limited to the shapes of the baseline hazards $h_{i0}(t)$ and the values of the log hazard ratios $\tilde{\beta}_i$. That being said, the overall treatment effect is not a readily defined value as a unique functional of all such baseline functions and parameters. For the first time in the literature, we point out that the problem has to be addressed in the misspecified model framework. Typically the target to be estimated with a misspecified model is the limit of a chosen estimator. Its definition has to depend on the modelling assumptions chosen by the researchers. We will demonstrate two kinds of such definitions in this paper and provide explaination for their very delicate differences. The second question naturally comes after the first one, i.e., how to generate a statistically efficient estimate for the overall drug efficacy after a proper definition? It is even more challenging if patient line data are not available but only aggregate statistics can be accessed for some of the trials.\\

\section{Combined hazard ratio as a limit of the MPLE from pooled data}
It is sufficient to consider only two independent trials. Generalization of the results to cover more trials is straightforward. Assume that $n$ patients were enrolled for the first clinical trial. The survival times and baseline demographics of the patients are denoted $\mathbf{X}=\{(X_i, \delta_i, \mathbf{Z}_i), i=1,\ldots,n\}$. Here $X_i$ is the right-censored survival time of the $i$-th patient, $\delta_i=0/1$ indicates that $X_i$ is censored or an observed event time, $\mathbf{Z}_i$ is a $k$-dimensional covariate with probability density function $f_\mathbf{Z}(z_1, \ldots, z_k)$. In many cases the distribution $f_\mathbf{Z}$ may vary by the trials. For simplicity, we assume the same distribution $f_\mathbf{Z}$ for all the trials in this paper. Extension of our methods to accommodate different distributions for the covariates is obvious. Assuming independent censoring, the patients' survival times are i.i.d.'s with proportional hazard
\[ h_X(x|\mathbf{Z})=h_{X0}(x)\exp(\tilde{\alpha}'\mathbf{Z}). \]
By definition, the pdf of an uncensored $X_i$ (i.e., $\delta_i=1$) conditioned on $\mathbf{Z}_i$ is
\[ h_{X0}(x)\exp(\tilde{\alpha}'\mathbf{Z})e^{-\exp(\tilde{\alpha}'\mathbf{Z})H_{X0}(x)}, \]
where $h_{X0}(\cdot)$ is an arbitrary baseline hazard function of the first trial and $H_{X0}(x)=\int_0^x h_{X0}(u)du$ is the corresponding culmulative hazard. Similar notations can be defined for the second trial. Let $\mathbf{Y}=\{(Y_j, \delta_j,\mathbf{Z}_j)$, $j=1,\ldots, m\}$ be the data collected from the second trial. Assuming proportional hazard and independent censoring, the patients from the second trial has hazards
\[ h_Y(y|\mathbf{Z})=h_{Y0}(y)\exp(\tilde{\beta}'\mathbf{Z}). \]
The trial-specific treatment effect estimates ($\hat{\alpha}$ and $\hat{\beta}$) and their asymptotic variance-covariance matrices can be easily estimated using the MPLEs calculated from the $\mathbf{X}$ and $\mathbf{Y}$ data respectively. We are particularly interested in the setup where the underlying true values of $\tilde{\alpha}\neq\tilde{\beta}$. Without loss of generality, let the first component of the covariates ($Z_{i1}$) be the arm indicator for the $i$-th patient. The patient is in the treatment arm if $Z_{i1}=1$ or he is in the control arm if $Z_{i1}=0$. The value of $\alpha_1$ and $\beta_1$, the first component of the regression parameter $\tilde{\alpha}$ and $\tilde{\beta}$ respectively, are of most concern as a measurement of the treatment effect. It is convenient to assume $\alpha_1<\beta_1$ for all the discussions presented in this paper.\\

\noindent
It is natural to consider the pooled patient line data  
\begin{eqnarray*} 
	&&\{(W_i,\delta_i,\mathbf{Z}_i), i=1,\ldots,n+m\} \\
	&=& \{(X_i,\delta_i,\mathbf{Z}_i), i=1,\ldots,n\}\cup\{(Y_j,\delta_j,\mathbf{Z}_j),j=1,\ldots,m\}.
\end{eqnarray*}
Most researchers consider the MPLE calculated from all $n+m$ (denoted by $N$) pooled patient line (\textit{PL}) records not only as an appropriate estimate for the overall treatment effects but also the best (especially when compared to the linear estimates $\hat{\beta}_L$ or $\hat{b}_L$) answer to our first question. The only obvious drawback of the MPLE is that it requires the knowledge of all $N$ patient line data:
\begin{eqnarray} 
	\hat{\theta}_{PL}=\arg\ \max_{\tilde{\theta}}\prod_{i=1}^{n+m}\left[\frac{\exp(\tilde{\theta}'\mathbf{Z}_i)}{\sum_{j\in\Re_i}\exp(\tilde{\theta}'\mathbf{Z}_j)}\right]^{\delta_i}, 
	\label{pl_log_hazard_ratio}
\end{eqnarray}
where $\Re_i$ is the set of labels for those patients (originally from $\mathbf{X}$ or $\mathbf{Y}$) who are at risk at time $W_i-$. \\


\noindent
The overall log hazard ratio $\tilde{\theta}$ is hence defined as the limit of $\hat{\theta}_{PL}$ when both $n$ and $m\rightarrow\infty$. It is worth to point out that one needs to first determine what is an appropriate statistic for the overall treatment effect based on both statstical and clinical thinking. Then the parameter to be estimated follows as the limit of the statistic, not vice versa. A different version of the overall treatment effect can be as valid based on other assumptions about the estimating procedure. We will extend the discussion to provide such an example in the Section 4.\\

\noindent
It is equivalent to imposing a misspecified Cox model (working model) on the pooled data such that the combined hazard can be written as
\begin{equation} 
	h_W(w|\mathbf{Z})=h_{W0}(w)\exp(\tilde{\theta}'\mathbf{Z}). \label{misspecified_model_hazard}
\end{equation}
Assume that the two trials have the same baseline hazard $h_{X0}(t)=h_{Y0}(t)=h_0(t)$ for all $t>0$ because the control arms are usually subject to standard of care. Such treatments are well established for the general population. The true pdf of the pooled data $W_i$ is a mixture of two proportional hazard models
\[ ph_0(w)e^{\tilde{\alpha}'\mathbf{Z}}e^{-\exp(\tilde{\alpha}'\mathbf{Z})H_0(w)} + (1-p)h_0(w)e^{\tilde{\beta}'\mathbf{Z}}e^{-\exp(\tilde{\beta}'\mathbf{Z})H_0(w)}, \]
where $n/(n+m)\rightarrow p$ as $n,m\rightarrow\infty$ is a fixed ratio of the sample sizes controlled by the researcher. The true model for $W_i$ does not satisfy the proportional hazard assumption:
\[  h(w|\mathbf{Z})=h_0(w)\frac{pe^{\tilde{\alpha}'\mathbf{Z}}e^{-\exp(\tilde{\alpha}'\mathbf{Z})H_0(w)}+(1-p)e^{\tilde{\beta}'\mathbf{Z}}e^{-\exp(\tilde{\beta}'\mathbf{Z})H_0(w)}}{pe^{-\exp(\tilde{\alpha}'\mathbf{Z})H_0(w)}+(1-p)e^{-\exp(\tilde{\beta}'\mathbf{Z})H_0(w)}} .\]
The formulation of the limit of $\hat{\theta}_{PL}$ can be studied using the techniques developed for misspecified models in Struthers and Kalbfleisch (1986) and Lin and Wei (1989). Let $h_i(t)$ be the true hazard function of the $i$-th patient from the pooled dataset $\mathbf{W}$ and $R_i(t)=\mathbbm{1}_{W_i\geq t}$ be the at-risk process at arbitrary time $t>0$, $i=1,\ldots, n+m$. It is convenient to define the notations following the convention of Andersen and Gill (1989), Struthers and Kalbfleisch (1986) and Lin and Wei (1989):
\begin{eqnarray*}
	s^{(0)}(t) = E\left[\sum_{i=1}^NR_i(t)h_i(t)\right], s^{(1)}(t) = E\left[\sum_{i=1}^NR_i(t)h_i(t)\mathbf{Z}_i\right],\\
	s^{(0)}(\tilde{\theta}, t) = E\left[\sum_{i=1}^N R_i(t)e^{\tilde{\theta}'\mathbf{Z}_i}\right], s^{(1)}(\tilde{\theta}, t) = E\left[\sum_{i=1}^N R_i(t)e^{\tilde{\theta}'\mathbf{Z}_i}\mathbf{Z}_i\right].
\end{eqnarray*}
Here the expected values are defined with respect to the true distribution of $W_i$ and $\mathbf{Z}_i$.\\

\noindent
\textit{\textbf{Proposition 1.} (based on \textbf{Theorem 2.1} of Lin and Wei (1989)) Let $\hat{\theta}_{PL}$ be the MPLE of the log hazard ratios for the overall treatment effect as defined in (\ref{pl_log_hazard_ratio}). When $n, m\rightarrow\infty$, $\hat{\theta}_{PL}$ converges in probability to the unique solution of the following equation:
\begin{equation} \int_0^\infty s^{(1)}(t)dt-\int_0^\infty\frac{s^{(1)}(\tilde{\theta},t)}{s^{(0)}(\tilde{\theta},t)}s^{(0)}(t)dt=0. \label{Lin_equation} \end{equation}}
$\Box$\\

\noindent 
Without censoring it follows the definition of $R_i(t)$ and the hazard $h_i(t)$ that
\[ E[R_i(t)h_i(t)|\mathbf{Z}_i]=P[W_i\geq t|\mathbf{Z}_i]h_i(t). \]
The right hand side of the above formula is simply the pdf of $W_i$ by the definition of the hazard $h_i(\cdot)$. It can be written as $f_X(t|\mathbf{Z})$ or $f_Y(t|\mathbf{Z})$ respectively,  in the form of
\[ h_{0}(t)\exp(\tilde{\alpha}'\mathbf{Z}_i)e^{-\exp(\tilde{\alpha}'\mathbf{Z}_i)H_{0}(t)}\ \mathrm{or}\ h_{0}(t)\exp(\tilde{\beta}'\mathbf{Z}_i)e^{-\exp(\tilde{\beta}'\mathbf{Z}_i)H_{0}(t)} \]
depending on whether the $i$-th patient is from the first or the second trial. Hence $s^{(0)}(t)$ and $s^{(1)}(t)$ can be simplified by the following calculation
\begin{eqnarray*}
	E[R_i(t)h_i(t)]&=&E[E[R_i(t)h_i(t)|\mathbf{Z}_i]]\\
	&=&E_{\mathbf{Z}}[f_X(t|\mathbf{Z})]\ \mathrm{if}\ W_i\in\mathbf{X},\ E_\mathbf{Z}[f_Y(t|\mathbf{Z})]\ \mathrm{otherwise}.
\end{eqnarray*}
Similarly $s^{(0)}(\tilde{\theta},t)$ and $s^{(1)}(\tilde{\theta},t)$ can be simplified using $E[R_i(t)|\mathbf{Z}_i]=P[W_i\geq t|\mathbf{Z}_i]$ in the form of 
\[ e^{-H_0(t)\exp(\tilde{\alpha}'\mathbf{Z}_i)}\ \mathrm{or}\ e^{-H_0(t)\exp(\tilde{\beta}'\mathbf{Z}_i)}.\]
Expand the shorthand notations defined for Proposition 1, we have
\begin{eqnarray*} 
	s^{(0)}(t) &=& E[nf_X(t|\mathbf{Z})+mf_Y(t|\mathbf{Z})],\\
 	s^{(1)}(t) &=& E[nf_X(t|\mathbf{Z})\mathbf{Z}+mf_Y(t|\mathbf{Z})\mathbf{Z}],\\
	s^{(0)}(\tilde{\theta}, t) &=& E[nP(X\geq t)e^{\tilde{\theta}'\mathbf{Z}}+mP(Y\geq t)e^{\tilde{\theta}'\mathbf{Z}}],\\
	s^{(1)}(\tilde{\theta}, t) &=& E[nP(X\geq t)e^{\tilde{\theta}'\mathbf{Z}}\mathbf{Z}+mP(Y\geq t)e^{\tilde{\theta}'\mathbf{Z}}\mathbf{Z}].
\end{eqnarray*}
All the above notations turn out to depend on no random variables other than the covariate $\mathbf{Z}_i$'s. The subscripts for $\mathbf{Z}$ are suppressed with the assumption that the $\mathbf{Z}_i$'s are independent and identically distributed. Under mild smoothness conditions, the first term of equation (\ref{Lin_equation}) can be simplified by switching the order of the integrals:
\[  \int_0^\infty s^{(1)}(t)dt=E\left[n\int_0^\infty f_X(t|\mathbf{Z})dt\mathbf{Z}+m\int_0^\infty f_Y(t|\mathbf{Z})dt\mathbf{Z}\right]=NE(\mathbf{Z}). \]
Plug in (\ref{Lin_equation}) with the definitions of the pdf's and survival functions. For sufficiently large $n$ and $m$, substitute $n/(n+m)$ by $p$, i.e., the fixed ratio of the study sizes. Eliminate $H_0(t)$ by letting $u=H_0(t)$ and hence $du=h_0(t)dt$. We derived an equation for the definition of the overall treatment effect.\\

\noindent
\textit{\textbf{Corollary 1.1} Assume no censoring in the data $\mathbf{X}$ and $\mathbf{Y}$. The true treatment effect (log hazard ratios) from either trials are known and denoted by $\tilde{\alpha}$ and $\tilde{\beta}$ respectively. The overall log hazard ratio $\theta^*_{PL}$ defined as the limit of $\hat{\theta}_{PL}$ with $n,m\rightarrow\infty$ is the unique solution to the following equation:
\begin{eqnarray}
 E(\mathbf{Z})=\int_0^\infty \frac{pE\left(e^{\tilde{\theta}'\mathbf{Z}-\exp(\tilde{\alpha}'\mathbf{Z})u}\mathbf{Z}\right) +(1-p)E\left(e^{\tilde{\theta}'\mathbf{Z}-\exp(\tilde{\beta}'\mathbf{Z})u}\mathbf{Z}\right)}{pE\left(e^{\tilde{\theta}'\mathbf{Z}-\exp(\tilde{\alpha}'\mathbf{Z})u}\right)+(1-p)E\left(e^{\tilde{\theta}'\mathbf{Z}-\exp(\tilde{\beta}'\mathbf{Z})u}\right)}\cdot\hspace{2cm}\label{Lin_solution}\\
 \hspace{2cm}\left[pE\left(e^{\tilde{\alpha}'\mathbf{Z}-\exp(\tilde{\alpha}'\mathbf{Z})u}\right)+(1-p)E\left(e^{\tilde{\beta}'\mathbf{Z}-\exp(\tilde{\beta}'\mathbf{Z})u}\right)\right]du. \nonumber
\end{eqnarray}}
$\Box$\\

\noindent
The distribution of the covariates $\mathbf{Z}$ can be well approximated using data from the general intent-to-treat population. Once the distribution of $\mathbf{Z}$ is known, equation (\ref{Lin_solution}) can be solved numerically.\\

\noindent
Usually one needs all $N$ patient line data to define the MPLE $\hat{\theta}_{PL}$ as well as its sandwich-type robust variance estimate. Nonetheless, if only aggregate data such as the MPLEs $\hat{\alpha}$ and $\hat{\beta}$ from individual trials are available for some reason, one may solve equation (\ref{Lin_solution}) with $\tilde{\alpha}$ and $\tilde{\beta}$ substituted by their estimates $\hat{\alpha}$ and $\hat{\beta}$. The solution to such an equation, which is henceforth denoted by $\hat{\theta}_{M}$ (the subscript $M$ stands for ``misspecified model''), does not rely on the knowledge of the baseline hazard $h_0(\cdot)$. It is a semiparametric, asymptotically efficient estimate to $\theta_{PL}^*$ because the MPLEs $\hat{\alpha}$ and $\hat{\beta}$ are based on maximum likelihoods. Such procedures are well known for being invariant with regard to functional transformations.\\

\noindent
The condition of no censoring in Corollary 1.1 is natural for the definition of $\theta_{PL}^*$ because we are only interested in the performance of the therapy.  Censoring is considered as noise imposed on the observed survival times and should be excluded from the estimating procedure if at all possible. The MPLEs $\hat{\alpha}$ and $\hat{\beta}$ reported for the individual trials are (asymtotically) unbiased for the underlying log hazard ratios $\tilde{\alpha}$ and $\tilde{\beta}$ even if the data $\mathbf{X}$ and $\mathbf{Y}$ are censored. Therefore $\hat{\theta}_M$ always remains an unbiased estimate for the overall treatment effect $\theta_{PL}^*$ no matter the data are censored or not. \\

\noindent
\textbf{Example 1.} Usually the effects of the covariates are assumed to be sorted out by proper randomization. One of the most important analysis used in practice use the treatment arm indicator $Z=0/1$ as the only covariate with $q=P(Z=1)$. Denote the hazard ratio of the first trial by $a=e^\alpha$, the hazard ratio of the second trial by $b=e^\beta$. By corollary 1.1, the MPLE $\hat{c}_{PL}=Exp(\hat{\theta}_{PL})$ calculated from the uncensored pooled line data of $N$ patient records converges to the solution of the following equation about $c$:
\begin{eqnarray} 
	1=\int_0^\infty \frac{(1-q)e^{-u}+pqae^{-au}+(1-p)qbe^{-bu}}{(1-q)e^{-u}+pqce^{-au}+(1-p)qce^{-bu}}\cdot\hspace{2cm} \label{Lin_solution_single_covariate}\\
	\hspace{3cm}\left[pce^{-au}+(1-p)ce^{-bu}\right]du. \nonumber
\end{eqnarray}
Equation (\ref{Lin_solution_single_covariate}) can be simplified after a series of basic algebraic transformations:
\begin{equation}
	1=\int_0^\infty \frac{(1-q)e^{-u}+pqae^{-au}+(1-p)qbe^{-bu}}{(1-q)e^{-u}+pqce^{-au}+(1-p)qce^{-bu}}\cdot e^{-u}du \label{Lin_solution_single_covariate_simple}
\end{equation}
It is easy to see that the right-hand-side of (\ref{Lin_solution_single_covariate_simple}) is a strictly decreasing, convex function of $c$. It implies that the overall hazard ratio defined by the solution $c^*_{PL}$ to equation (\ref{Lin_solution_single_covariate_simple}) must take a value on the open interval $(a,b)$. In Appendix 5.2, we will also show that in most cases $c^*_{PL}$ is superior to (i.e., smaller than) the linear alternative $\exp[p\alpha+(1-p)\beta]$ and $pa+(1-p)b$. $\Box$\\

\noindent
Now take a closer look at $\hat{\theta}_{PL}$. It is conventionally considered as the best estimate to the overall log hazard ratio of the combined trials. However, $\theta_{PL}^*$ is actually defined to be the limit of the uncensored $\hat{\theta}_{PL}$. The estimator $\hat{\theta}_{PL}$ is subject to the impact of censoring and can be biased for $\theta_{PL}^*$. By (\ref{Lin_solution}), the value of $\theta^*_{PL}$ only relies on $\tilde{\alpha}$, $\tilde{\beta}$ and the distributions of $\mathbf{Z}$.\\

\noindent
\textit{\textbf{Corollary 1.2} Assume independent right censoring for both $\mathbf{X}$ and $\mathbf{Y}$ such that the survival function of the censoring times are denoted $C_X(t|\mathbf{Z})$ and $C_Y(t|\mathbf{Z})$ respectively. The limit of $\hat{\theta}_{PL}$ with $n,m\rightarrow\infty$ is the unique solution to the following equation with respect to $\tilde{\theta}$ with known values of $\tilde{\alpha}$ and $\tilde{\beta}$:
\begin{equation}
 \int_0^\infty E[pf_X(t|\mathbf{Z})C_X(t|\mathbf{Z})\mathbf{Z}+(1-p)f_Y(t|\mathbf{Z})C_Y(t|\mathbf{Z})\mathbf{Z}] dt= \label{corollary_censor}
\end{equation}
\begin{eqnarray}
 \int_0^\infty \frac{pE\left(e^{\tilde{\theta}'\mathbf{Z}-\exp(\tilde{\alpha}'\mathbf{Z})H_0(t)}\mathbf{Z}C_X(t|\mathbf{Z})\right) +(1-p)E\left(e^{\tilde{\theta}'\mathbf{Z}-\exp(\tilde{\beta}'\mathbf{Z})H_0(t)}\mathbf{Z}C_Y(t|\mathbf{Z})\right)}{pE\left(e^{\tilde{\theta}'\mathbf{Z}-\exp(\tilde{\alpha}'\mathbf{Z})H_0(t)}C_X(t|\mathbf{Z})\right)+(1-p)E\left(e^{\tilde{\theta}'\mathbf{Z}-\exp(\tilde{\beta}'\mathbf{Z})H_0(t)}C_Y(t|\mathbf{Z})\right)}\cdot\hspace{2cm}\nonumber\\
 \hspace{2cm}\left[pE\left(f_X(t|\mathbf{Z})C_X(t|\mathbf{Z})\right)+(1-p)E\left(f_Y(t|\mathbf{Z})C_Y(t|\mathbf{Z})\right)\right]dt.\hspace{2cm}. \nonumber
\end{eqnarray}}
$\Box$\\

\noindent
With a mixed population, censoring contains information about the source of the data, which is correlated with the length of the subject's expected survival time beyond censoring. The limit of $\hat{\theta}_{PL}$ must contain all such information as the censoring mechanism $C_X(t|\mathbf{Z})$, $C_Y(t|\mathbf{Z})$ and the baseline hazard $H_0(t)$. \\

\noindent
\textbf{Example 2.} Assume the same setup as in Example 1. The survival times observed from either trials follow Cox models defined with a treatment arm indicator $Z\sim Bernoulli(q)$. Both clinical trials will be terminated at given time $T_{max}>0$. Hence the censoring time has a point mass of unity at $T_{max}$. That is, $C_X(t|\mathbf{Z})=C_Y(t|\mathbf{Z})=1$ if $t<T_{max}$ and 0 otherwise. Expand the definitions of the distribution functions and use the change-of-variable technique by letting $u=H_0(t)$ in equation (\ref{corollary_censor}), it becomes
\begin{equation}
	1-e^{-H_0(T_{max})}=\int^{H_0(T_{max})}_0 \frac{(1-q)e^{-u}+pqae^{-au}+(1-p)qbe^{-bu}}{(1-q)e^{-u}+pace^{-au}+(1-p)qce^{-bu}}\cdot e^{-u}du. \label{equation_censor}
\end{equation}
The limit of the MPLE $\hat{\theta}_{PL}$ calculated from the censored line data is the solution $c^*$ to the above equation. To study the bias in $\hat{\theta}_{PL}$, we need to compare $c^*$ against $c^*_{PL}$, which is the solution to equation (\ref{Lin_solution_single_covariate_simple}). For simplicity, denote the integrand in equation (\ref{equation_censor}) by $g(u|c)$ for fixed $a$ and $b$. Equation (\ref{Lin_solution_single_covariate_simple}) asserts that $g(u|c^*_{PL})$ is a well defined pdf because its integral on $(0, \infty)$ equals unity. Consider the fact that $c^*_{PL} < pa+(1-p)b$ (Appendix 5.2). It is easy to prove that the fraction term
\[  \frac{(1-q)e^{-u}+pqae^{-au}+(1-p)qbe^{-bu}}{(1-q)e^{-u}+pac^*_{PL}e^{-au}+(1-p)qc^*_{PL}e^{-bu}} >1 \]
if $0<u < [\ln(1-p)(c^*_{PL}-b)/p/(a-c^*_{PL})]/(b-a)$ and it is less than 1 otherwise. Hence the pdf $g(u|c^*_{PL})$ only intersects a standard Exponential pdf $e^{-u}$ at one point, which in turn implies that the distribution defined by $g(u|c^*_{PL})$ is stochastically smaller than a standard Exponential distribution (Fill and Machida 2001). Therefore, the cdf corresponding to $g(u|c^*_{PL})$ is always bigger than that of the standard Exponential distribution $\int_0^H e^{-u}du=1-e^{-H}$ for any given $H$:
\[ 1-e^{-H_0(T_{max})}<\int^{H_0(T_{max})}_0 \frac{(1-q)e^{-u}+pqae^{-au}+(1-p)qbe^{-bu}}{(1-q)e^{-u}+pac^*_{PL}e^{-au}+(1-p)qc^*_{PL}e^{-bu}}\cdot e^{-u}du. \]
Note that $g(u|c)$ is monotonically decreasing with respect to $c$. To make an equality as in equation (\ref{equation_censor}), its solution $c^*$ must be greater than $c^*_{PL}$. When the survival time data are censored at a maximum allowable trial length $T_{max}$, $\hat{\theta}_{PL}$ is always associated with a positive bias. The bias decreases with $T_{max}$. $\Box$\\

\noindent
Corollary 1.2 indicates that the most accepted MPLE $\hat{\theta}_{PL}$ is not robust against variability in the treatment effects observed from multiple trials. Hence we recommend reporting the overall log hazard ratio for multiple trials using $\hat{\theta}_M$ rather than $\hat{\theta}_{PL}$. The former is not only robust (unbiased despite of censoring) but also provides better chance for the researchers because it only requires aggregate statistics from each sub-population of concern. The following example illustrates the impact of censoring on $\hat{\theta}_{PL}$ using simulated data.\\

\noindent
\textbf{Example 3.} We performed 1000 rounds of independent simulations to mimic the following scenario: the survival times of 200 patients treated in the first trial follow an $Exp(0.3)$ distribution. With 1:1 randomization (i.e., $q=0.5$), another 200 patients in the control group has survival times sampled from a standard exponential distribution $Exp(1)$. Let $p=0.7$, the second trial enrolls 85 patients for the treated group and 85 patients for the control group. To mimic a lower drug efficacy, the survival times of the treated patients from the second trial are sampled from an $Exp(0.8)$ distribution, while the survival times of the control group patients are sampled from $Exp(1)$.\\

\noindent
Without censoring, we can calculate $\hat{\theta}_{PL}$ for each of the 1000 simulated data set and summarize the distribution of $\hat{\theta}_{PL}$ using its empirical distribution. In this example, it was reported that $E(\hat{\theta}_{PL})=-0.926$ (equivalent to $\ln(0.396)$) with a 95\% confidence interval of $(-1.088, -0.756)$. According to Corollary 1.1 this is an (asymtotically) unbiased estimate for the true overall log hazard ratio $\theta^*_{PL}$ as there is no censoring.\\

\noindent
We continued to censor the 1000 simulated data sets using various censoring time $T_{max} \in [1,10]$ and tried to calculate the estimates $\hat{\theta}_M$ and $\hat{\theta}_{PL}$ respectively for each specific $T_{max}$. In Figure 1(i) $\hat{\theta}_M$'s are reported as the solution to equation (\ref{Lin_solution_single_covariate_simple}) with $\alpha$ and $\beta$ substibuted by $\hat{\alpha}$ and $\hat{\beta}$, which are calculated from the censored trial 1 and trial 2 data respectively. On the other hand, in Figure 1(ii) the MPLE $\hat{\theta}_{PL}$'s are calculated using (\ref{pl_log_hazard_ratio}) with all 570 censored line data from either trials pooled together. \\
\begin{figure}[ht!]
\centering
\includegraphics[width=130mm]{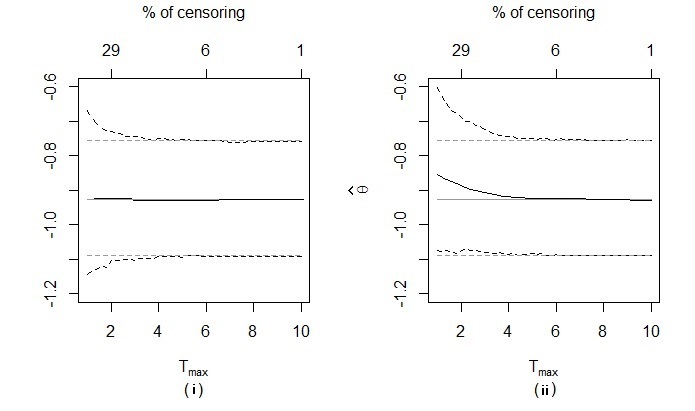}
\caption{(i) Mean and the 95\% CI for $\hat{\theta}_M$. (ii) Mean and the 95\% CI for $\hat{\theta}_{PL}$.}
\end{figure}

\noindent
The gray lines in Figure 1 outline the true overall hazard ratio (-0.926) and its 95\% confidence interval $(-1.088, -0.756)$. The bias of $\hat{\theta}_{PL}$ can be easily observed in Figure 1(ii). When the trials are censored at $T_{max}=1$, about 51\% of the collected data are censored. The MPLE $\hat{\theta}_{PL}$ reported for this scenario has a mean of -0.854 and a 95\% confidence intervals of $(-1.074, -0.601)$. It accounts for about 7.8\% of positive bias in the log hazard ratio estimates. Such a result can be of serious concern for those clinical trials expected to be associated with low event rates, e.g., a trial for breast cancer treatments. For the same censored dataset, $\hat{\theta}_M$ displayed in Figure 1(i) appears to be unbiased for $\theta^*_{PL}$. The confidence interval of $\hat{\theta}_M$ tends to be wider with more data being censored because the variance of the log hazard ratio estimate is proportional to the inverse of the number of observed events (Kalbfleisch and Prentice (2002)). $\Box$

\section{Combined hazard ratio via harmonic means}
We noted that $\hat{\theta}_{PL}$ and $\hat{\theta}_M$ are semi-parametric estimators without requiring any knowlege about the baseline. However, the baseline hazard $h_{W0}(\cdot)$ of the misspecified model (\ref{misspecified_model_hazard}) is different from the baseline hazard $h_0(\cdot)$ defined for each individual clinical trial. It is due to the fact that the MPLE procedure is the result of simultaneous maximization of the unknown parameter $\tilde{\theta}$ and the discretized baseline hazard $h_{W0}(\cdot)$ in the form of a Breslow hazard estimate (Breslow (1972), Johansen (1983)). The shape of the baseline hazard $h_0(\cdot)$ is twisted to fit the mixed event times and in return leads to a slightly under-estimated therapy effect, or equivalently an over-estimated log hazard ratio $\tilde{\theta}^*_{PL}$ in many cases. We will propose another method to avoid such undesirable effect. \\

\noindent
\textbf{Example 4.} Let $X_1, \ldots, X_n\sim Exp(a)$ and $X_{n+1}, \ldots, X_{2n}\sim Exp(1)$ be the i.i.d. event times observed from the treatment and control groups of the first trial,  $Y_1, \ldots, Y_m\sim Exp(b)$ and $Y_{m+1}, \ldots, Y_{2m}\sim Exp(1)$ be the event times recorded from the treatment and control groups of the second trial. Without loss of generality, assume $a<b<1$. No censoring is allowed for simplicity. Provided the complete patient line data, it is easy to estimate the treatment effects of the therapy in either clinical trial using the MLEs $\hat{a}=n/\sum_{i=1}^n X_i$ and $\hat{b}=m/\sum_{j=1}^m Y_j$. Assuming a misspecified $Exp(c)$ model for the pooled data, it is natural to estimate the overall treatment effect using 
\begin{eqnarray*}
	\hat{c}=\frac{n+m}{\sum_{i=1}^nX_i+\sum_{j=1}^mY_m}=\frac{1}{n/(n+m)\hat{a}+m/(n+m)\hat{b}}\\
	\stackrel{p}{\longrightarrow} \frac{1}{p/a+(1-p)/b}.
\end{eqnarray*}
Again, the overall hazard ratio can be defined as the limit of the chosen estimate $\hat{c}$. It turns out to be the harmonic mean of the individual trial effects $a$ and $b$ weighted by the ratio of the study sizes $p:(1-p)$. As a matter of fact, the harmonic mean effect is smaller than that defined by the semiparametric MPLE method in most cases (Appendix 5.2):
\[ \frac{1}{p/a+(1-p)/b}< c^*_{PL},\ \mathrm{if}\ a< b<1, \]
where $c^*_{PL}$ is the solution of equation (\ref{Lin_solution_single_covariate_simple}) as well as the limit of the MPLE $\exp(\hat{\theta}_{PL})$ calculated using $\{W_1,\ldots,W_{2n+2m}\}=\{X_1\ldots, X_{2n}, Y_1, \ldots, Y_{2m}\}$ without assuming an underlying exponential hazard. The difference between the two versions of the overall treatment effect estimate can be attributed to a twisted baseline hazard $h_{W0}(\cdot)$ approximate due to the MPLE procedure. Consider the two trials respectively, the Breslow hazard estimates 
\[ \hat{H}_0(t)= \sum_{i:X_i\leq t}\frac{1}{\sum_{\ell=1}^{2n}R_{X\ell}(X_i)\hat{a}^{Z_\ell}}\ \mathrm{or}\ \sum_{j:Y_j\leq t}\frac{1}{\sum_{\ell=1}^{2m}R_{Y\ell}(Y_j)\hat{b}^{Z_\ell}}\]
calculated from either sets of patient line data are both asymptotically unbiased for the true underlying linear culmulative hazard function $H_0(t)=t$. Here the only covariate $Z_\ell=0$ if the patient is in the control group and 1 if the patient is in the treatment group. The numerator of the fraction is always equal to one because there are no censoring. However, the Breslow estimate calculated from the pooled data $W_i$'s does not lead to a constant hazard estimate over time. At arbitrary time $t$, the non-parametric baseline hazard $h_{W0}(t)$ as a limit of the Breslow estimate is
\begin{equation} 
	h_{W0}(t)=\frac{pe^{-at}a+(1-p)e^{-bt}b+e^{-t}}{[pe^{-at}+(1-p)e^{-bt}]c^*_{PL}+e^{-t}}.
	\label{combined_Breslow_limit}
\end{equation}
The derivation of the above formula is provided in Appendix 5.1. Note that
\[
	h_{W0}(0)=\frac{pa+(1-p)b+1}{c^*_{PL}+1}\ \mathrm{and}\ h_{W0}(\infty)=\frac{a}{c^*_{PL}}.
\]
Unless $a=b$, the shape of $h_{W0}(t)$ is tilted to the right by the fact that $pa+(1-p)b>c^*_{PL}$ and $a<c^*_{PL}$ (Appendix 5.2). For sufficiently large $t$, $h_{W0}(t)$ is monotonically decreasing and $h_{W0}(t)=1$ happens at only one time point $t$. The trend in $h_{W0}(t)$ is consistent with the monotone changes in the mixture proportion of the complete population. At the very beginning, the estimated baseline hazard is bigger than one. The proportion of treated patients from the first trial increases with time because the patients from the second trial have shorter expected life ($a<b$) and finally the baseline hazard estimate is dominated by $a$, the effect of the first trial, and scaled by $1/c^*_{PL}$ such that the average baseline hazard is close to one. Hence a smaller combined effect (i.e., bigger hazard ratio $c^*_{PL}$) is given by the MPLE procedure compared to the harmonic mean effect based on the parametric Exponential model. In a sense, the MPLE procedure is an overfit to the data if the researcher is confident about the fact that the control groups from either trial are not essentially different. $\Box$\\

\noindent
It is curious to see that the harmonic mean type of definition for the combined trial effect can be extended to address statistical models assuming much more general conditions where neither the underlying exponential distribution nor the univariate covariate structure is needed. Let $\mathbf{X}$ denote the observed survival times of $n$ patients recruited for a clinical trial. Their uncensored survival times follow a proportional hazard model with hazard $h_X(t|\mathbf{Z})=h_0(t)\exp(\tilde{\alpha}'\mathbf{Z})$. Similarly, the survival times ($\mathbf{Y}$) of another $m=n(1-p)/p$ patients recruited for a second trial also follow a proportional hazard model with the same baseline hazard $h_0(t)$ and a log hazard ratio $\tilde{\beta}$. Here $\tilde{\alpha}$ is a $k$-variate vector of the same structure as $\tilde{\beta}$ but of different values. The formulation of the baseline hazard $h_0(t)$ is unknown. The MPLE is not appropriate for the estimation of the log hazard ratio ($\tilde{\theta}$) for the combined population if one needs to avoid a twisted baseline hazard. Instead, we set out to define the MLE for $\tilde{\theta}$. It leads to a different version of a semi-parametric estimate ($\hat{\theta}_{HM}$, where the subscript $HM$ stands for ``harmonic mean'') for the overall log hazard ratio. It is based on aggregate statistics only. Patient line data are not required for the realization of $\hat{\theta}_{HM}$.\\

\noindent
\textit{\textbf{Proposition 2.} 
 Let $\hat{\theta}_{MLE}$ be the maximum likelihood estimate (MLE) for the overall log hazard ratio $\tilde{\theta}$ when a proportional hazard model with the following pdf is fitted to the pooled $N=n+m$ surival time records:
\begin{equation}
	f(t;\tilde{\theta}|\mathbf{Z})=h_0(t)e^{\tilde{\theta}'\mathbf{Z}}e^{-\exp(\tilde{\theta}'\mathbf{Z})H_0(t)}, t\geq 0.\label{misspecified_model}
\end{equation}
When $n,m\rightarrow\infty$, $\hat{\theta}_{MLE}$ converges in probability to a constant $\theta^*_{HM}$, which is the unique solution to the following equation with respect to $\tilde{\theta}$:
\begin{equation}
	E(\mathbf{Z})=E\left[e^{\tilde{\theta}'\mathbf{Z}}\mathbf{Z}\left(\frac{p}{e^{\tilde{\alpha}'\mathbf{Z}}}+\frac{1-p}{e^{\tilde{\beta}'\mathbf{Z}}}\right)\right]. \label{hm_equation}
\end{equation}
}

\noindent
\textit{Proof.} Again, we assume no censoring when trying to define the ``true'' values of the overall treatment effect because no information other than the measurements of the treatment effect should be of concern. Such conditions can be loosen when we get to the discussions about the estimating procedures for the treatment effects established here.\\

\noindent 
Assuming the misspecified proportional hazard model (\ref{misspecified_model}), the joint pdf of the $N$ observed survival times is
\[ \prod_{i=1}^Nh_0(T_i)e^{\tilde{\theta}'\mathbf{Z}_i}e^{-e^{\tilde{\theta}'\mathbf{Z}_i}H_0(T_i)}. \]
The derivative (w.r.t. $\tilde{\theta}$) of the logarithm of the misspecified joint pdf is 
\[ \nabla_{\tilde{\theta}}\ell(T_1,\ldots,T_N,\mathbf{Z}_1,\ldots,\mathbf{Z}_N)= \sum_{i=1}^N \mathbf{Z}_i-e^{\tilde{\theta}'\mathbf{Z}_i}\mathbf{Z}_iH_0(T_i). \]
By the law of large number, 
\[ \nabla_{\tilde{\theta}}\ell(T_1,\ldots,T_N,\mathbf{Z}_1,\ldots,\mathbf{Z}_N)/N \stackrel{p}{\longrightarrow} E_{mix}\left[\mathbf{Z}-e^{\tilde{\theta}'\mathbf{Z}}\mathbf{Z}H_0(T)\right], \]
where $E_{mix}(\cdot)$ stands for the expectation defined in terms of the true mixed model with pdf $f_{mix}$:
\[ f_{mix}(t,\mathbf{z})= [ph_0(t)e^{\tilde{\alpha}'\mathbf{z}}e^{-e^{\tilde{\alpha}'\mathbf{z}}H_0(t)}+(1-p)h_0(t)e^{\tilde{\beta}'\mathbf{z}}e^{-e^{\tilde{\beta}'\mathbf{z}}H_0(t)}]f_\mathbf{Z}(\mathbf{z}).\]
This is a typical setup of an estimating equation. Under mild regularity conditions, it can be proved that the MLE $\hat{\theta}$ defined by the solution to the equation $\nabla_{\tilde{\theta}}\ell=0$ converges in probability to a $\tilde{\theta}_{HM}$, which is the solution to 
\begin{equation}
	E_{mix}\left[\mathbf{Z}-e^{\tilde{\theta}'\mathbf{Z}}\mathbf{Z}H_0(T)\right]=0. \label{expectation_equation_hm}
\end{equation}
Detail discussions about the asymptotics of the MLE derived from misspecified models are available in White (1982). Note that (\ref{expectation_equation_hm}) is equivalent to
\begin{eqnarray} 
	E_{mix}(\mathbf{Z})&=&E_{mix}\left[e^{\tilde{\theta}'\mathbf{Z}}\mathbf{Z}H_0(T)\right]\nonumber\\
	&=& E_{\mathbf{Z}}\left[E_{mix}\left(e^{\tilde{\theta}'\mathbf{Z}}\mathbf{Z}H_0(T)|\mathbf{Z}\right)\right]\nonumber\\
	&=& E_{\mathbf{Z}}\left[e^{\tilde{\theta}'\mathbf{Z}}\mathbf{Z}E_{mix}(H_0(T)|\mathbf{Z})\right]. \label{proof_corollary2}
\end{eqnarray}
To simplify the right hand side of (\ref{proof_corollary2}), we calculate
\begin{eqnarray*}
	&& E_{mix}[H_0(T)|\mathbf{Z}]\\
	&=&\int_0^\infty H_0(t)f_{mix}(t|\mathbf{Z})dt\\
	&=&\int_0^\infty H_0(t)\left[ph_0(t)e^{\tilde{\alpha}'\mathbf{Z}}e^{-e^{\tilde{\alpha}'\mathbf{Z}}H_0(t)}+(1-p)h_0(t)e^{\tilde{\beta}'\mathbf{Z}}e^{-e^{\tilde{\beta}'\mathbf{Z}}H_0(t)}\right]dt\\
	&=&\int_0^\infty u\left[pae^{-au}+(1-p)be^{-bu}\right]du\\
	&=&p/a+(1-p)/b.
\end{eqnarray*}
Note that $h_0(t)dt=dH_0(t)$. The third line of the above equation was due to the following definition of notations:  $H_0(t)=u$, $e^{\tilde{\alpha}'\mathbf{Z}}=a$ and $e^{\tilde{\beta}'\mathbf{Z}}=b$. Equation (\ref{hm_equation}) follows (\ref{proof_corollary2}) with the definition of $E_{mix}[H_0(T)|\mathbf{Z}]$ plugged in. $\Box$\\

\noindent
Following the discussions about misspecified models in White (1982) and Akaike (1973), it can be seen that equation (\ref{expectation_equation_hm}) actually defines $\theta^*_{HM}$ as the maximizer of the following expectation:
\[ E_{mix}[\ln(f(T;\tilde{\theta}|\mathbf{Z})f_\mathbf{Z}(\mathbf{Z}))]. \]
Hence, $\theta^*_{HM}$ has an obvious geometric interpretation. It minimizes the Kullback-Leibler distance between the candidate working models $f(t;\tilde{\theta}|\mathbf{Z})$ and the true model $f_{mix}(t|\mathbf{Z})$:
\begin{eqnarray*}
	\theta^*_{HM}&=&\arg \min_{\tilde{\theta}} D[f_{mix}, f(t;\tilde{\theta})]\\
	&=& \arg \min_{\tilde{\theta}} E_{mix}[\ln(f_{mix}(T,\mathbf{Z}))]-E_{mix}[\ln(f(T,\tilde{\theta},\mathbf{Z}))]. 
\end{eqnarray*}
Now we have another version of the definition for the overall log hazard ratio $\theta^*_{HM}$. It can be slightly smaller than $\theta^*_{PL}$ in most cases (Appendix 5.2). Both of these two numbers are valid measurements of the treatment effect on the combined population, though they are based on different modelling assumptions. We consider $\theta^*_{HM}$ a number closely related to the harmonic mean. Such an idea can be illustated by the following example.\\

\noindent
\textbf{Example 5.} Assume the same setup as in Example 1. The treatment group indicator $Z\sim Bernoulli(q)$ is the only covariate collected for the enrolled patients. The observed survival times (with or without censoring) follow proportional hazard models with hazard ratios $a=e^\alpha$ and $b=e^\beta$ respectively for trial 1 and 2. By equation (\ref{hm_equation}), the overall hazard ratio $c^*_{HM}$ is the solution to the following equation:
\[ q=E(Z)=E\left[c^ZZ\left(\frac{p}{a^Z}+\frac{1-p}{b^Z}\right)\right]=qc\left(\frac{p}{a}+\frac{1-p}{b}\right) \]
\begin{equation}
	\Longrightarrow c^*_{HM}=\frac{1}{p/a+(1-p)/b}. \label{c_hm}
\end{equation}
$\Box$\\

\noindent
Note that exponential survival times are not required in Example 5. The harmonic type of calculation (\ref{c_hm}) is applicable to any general proportional hazard modelling setup.\\

\noindent
Substitute $\alpha$ and $\beta$ by the corresponding MPLEs $\hat{\alpha}$ and $\hat{\beta}$ in equation (\ref{hm_equation}). Solve the estimating equation about $\tilde{\theta}$ and denote the solution by $\hat{\theta}_{HM}$.  It is an asymtotically efficient estimate for $\theta^*_{HM}$ because both $\hat{\alpha}$ and $\hat{\beta}$ are based on maximum likelihood principles and hence are invariant to any functional transformations. The asymptotic variance of $\hat{\theta}_{HM}$ can be derived using the delta method. Consider $\hat{\theta}_{HM}$ as an implicit function of $\hat{\alpha}$ and $\hat{\beta}$ by the estimating equation (\ref{hm_equation}). Let $\partial_j\tilde{\theta}=(\partial\theta_1/\partial\alpha_j,\ldots, \partial\theta_k/\partial\alpha_j)$, $\forall j=1, \ldots, k$. Assuming mild regularity conditions, e.g., dominated convergence for the variables defined in (\ref{hm_equation}), one may switch the order of integrations and differentiations. The values of $\partial_j\tilde{\theta}$ can be calculated by solving the following linear system about $\partial_j\tilde{\theta}$:
\[ E\left[e^{\tilde{\theta}'\mathbf{Z}}\mathbf{Z}\left(\frac{p}{e^{\tilde{\alpha}'\mathbf{Z}}}+\frac{1-p}{e^{\tilde{\beta}'\mathbf{Z}}}\right)\left(\mathbf{Z}'\partial_j\tilde{\theta}\right)\right]=E\left[ e^{\tilde{\theta}'\mathbf{Z}}\mathbf{Z}\frac{pZ_j}{e^{\tilde{\alpha}'\mathbf{Z}}}\right]. \]
Similarly it is easy to derive the formulas for $(\partial\theta_1/\partial\beta_j, \ldots, \partial\theta_k/\partial\beta_j)$, $\forall j=1,\ldots,k$. Usually the variance-covariance matrices of $\hat{\alpha}$ and $\hat{\beta}$, denoted by $Var(\hat{\alpha}$ and $Var(\hat{\beta})$, are reported together with the point estimate values. By the delta method, the asymtotic variance of $\hat{\theta}_{HM}$ can be calculated:
\[ \left( \begin{array}{ccc}
\partial\theta_1/\partial\alpha_1 & \ldots & \partial\theta_1/\partial\beta_k \\
 & \ldots &  \\
\partial\theta_k/\partial\alpha_1 & \ldots & \partial\theta_k/\partial\beta_k \end{array} \right)
\left( \begin{array}{cc}
Var(\hat{\alpha}) & \mathbf{0}_{k \times k} \\
\mathbf{0}_{k\times k} & Var(\hat{\beta}) \end{array} \right)
\left( \begin{array}{ccc}
\partial\theta_1/\partial\alpha_1 & \ldots & \partial\theta_k/\partial\alpha_1 \\
 & \ldots &  \\
\partial\theta_1/\partial\beta_k & \ldots & \partial\theta_k/\partial\beta_k \end{array} \right).
\]
A wald test can be developed for the values of $\theta^*_{HM}$ with its variance-covariance matrix calculated as above.\\

\noindent
\textbf{Example 5. (cont.)} When there is only one Bernoulli covariate $Z$ as was specified in Example 5, the variance of the log hazard ratio estimate
\[ \hat{\theta}_{HM}=\ln\frac{1}{p/\exp(\hat{\alpha})+(1-p)/\exp(\hat{\beta})} \]
can be calculated using the delta method:
\[ Var(\hat{\theta}_{HM})=\frac{p^2e^{-2\hat{\alpha}}Var(\hat{\alpha})+(1-p)^2e^{-2\hat{\beta}}Var(\hat{\beta})}{(pe^{-\hat{\alpha}}+(1-p)e^{-\hat{\beta}})^2}. \]
When $\alpha=\beta$, the above formula degenerages to $Var(\hat{\theta}_{HM})\approx p^2Var(\hat{\alpha})+(1-p)^2Var(\hat{\beta})$, implying that $\hat{\theta}_{HM}$ is equivalent to $\hat{\theta}_L=p\hat{\alpha}+(1-p)\hat{\beta}$ in this case.\\

\noindent
A Wald test for the overall treatment effect $H_0:\theta^*_{HM}=0$ can be defined  using the test statistic $\hat{\theta}/\sqrt{Var(\hat{\theta}_{HM})}$. $\Box$

\section{Conclusion}
In this paper we investigated various methods for the estimation of the overall treatment effect observed from a mixed patient population. Linear estimators in the form of $\hat{\theta}_L=\sum_i w_i\hat{\beta}_i$ or $\hat{c}_L=\sum_i w_i e^{\hat{\beta}_i}$ have been the favorite of many researchers for their simplicity. However, it is not mathematically justifiable to approximate the notoriously nonlinear hazard ratio using any linear estimators if the patient responses to the treatment are highly diversified in various sub-groups of the intent-to-treat population. In particular, we showed that $\hat{c}_L>\exp(\hat{\theta}_L)$ and both of them are, in most cases, positively biased for the hazard of the combined treated patients.\\

\noindent
We propose that an appropriate definition of the overall treatment effect for a mixed population should be first of all based on an estimating procedure that is justifiable from either a clinical or statistical perspective. The first candidate meets such criterion is the MPLE $\hat{\theta}_{PL}$ calculated from the pooled patient line data. It converges to a well-defined overall log hazard ratio $\theta^*_{PL}$ for the combined trials if the observed event times are not censored. However, $\hat{\theta}_{PL}$ is biased if the data are censored as in most of the real life examples.\\

\noindent
The MPLE $\hat{\theta}_{PL}$ has a robust version $\hat{\theta}_M$ defined with a misspecified proportional hazard model. It is an asymptotically efficient semi-parametric estimator calculated from the aggregate statistics $\hat{\alpha}$ and $\hat{\beta}$. It converges to $\theta^*_{PL}$ with increasing sample sizes despite censoring in the data. We noted that the baseline hazard function is twisted when applying the MPLE procedure to the observed survival times. To avoid tampering the shape of the non-parametric baseline hazard, we proposed a harmonic mean type of estimator $\hat{\theta}_{HM}$. Again, it is a semi-parametric estimator based on $\hat{\alpha}$ and $\hat{\beta}$ only. It converges to a $\theta^*_{HM}$, which minimizes the Kullback-Leibler distance between the true mixed proportional hazard model and the misspecified working model. The variance-covariance matrix of $\hat{\theta}_{HM}$ can be calculated using the delta method and the reported $Var(\hat{\alpha})$ and $Var(\hat{\beta})$. We also derived a Wald test for the values of $\theta^*_{HM}$.

\section{Appendix}
\subsection{Example 4 (cont.)}
In Example 4 of section 3, we noted that the MPLE procedure together with the Breslow estimate can twist the shape of the non-parametric baseline hazard function. Details of the calculations are provided here.\\

\noindent
Assume the same Exponential setup delineated in Example 4. At arbitrary time $t>0$, let $w$ be the event time (it can be observed from either trial) right before  $t$ and $w_+$ be the next event time. The Breslow hazard estimate for the infinitestimal time interval $(w, w_+)$ is
\begin{equation} 
	d\hat{\Lambda}_0(t)=\frac{1}{[n_{X}(t)+m_Y(t)]\hat{c}_{PL}+(n_{Xc}(t)+m_{Yc}(t))}, \label{combined_Breslow}
\end{equation}
where $n_X(t)$ and $m_Y(t)$ are the number of treated patients still at risk up to time $t$, $n_{Xc}$ and $m_{Yc}$ are the number of living control group patients. The formulation of $d\hat{\Lambda}_0(t)$ is discussed in many research papers and texts, e.g., Breslow (1972), Kalbfleisch and Prentice (2002) and Hanley (2008), though the meaning of the formulas remains unclear to many statisticians. It is corresponding to a discretized Poisson process with constant hazard between consecutive event times. All history up to time $w$ can be ignored since the hazard is defined as a conditional probability for the future beyond $w$. Consider a control group patient being alive at time $w$. The estimate $d\hat{\Lambda}_0(t)$ is in fact an empirical approximate for the probability of observing this control group patient die within the time interval $[w, w_+)$: $d\hat{\Lambda}_0(t)= P[X_c\leq (w_+-w)]$, where $X_c\sim Exp(\lambda_0)$ is the event time of the imaginary control group patient with a to-be-estimated intensity parameter $\lambda_0$ specifically defined for the time interval $[w, w_+)$. This is consistent with the definition of the cumulative hazard:
\[ \int_{w}^{w_+}\lambda_0(t)dt=\ln\frac{S_0(w_+)}{S_0(w)}\approx \frac{S_0(w_+)-S_0(w)}{S_0(w)}=P[X_c\leq (w_+-w)], \]
where $S_0(\cdot)$ denotes the baseline survival function. The length of the time interval $(w_+-w)$ is also exponentially distributed, according to the true underlying distribution, with an intensity of 
\[ n_X(t)a+m_Y(t)b+(n_{Xc}(t)+m_{Yc}(t)) \]
because each of the $n_X(t)$ living treated patients from trial 1, the $m_Y(t)$ treated patients from trial 2 and the $(n_{Xc}(t)+m_{Yc}(t))$ control group patients can be considered as a competing risk. Use the joint pdf of two independent exponential random variables ($X_c$ vs. $w_+-w$) to calculate that 
\begin{eqnarray*}
	d\hat{\Lambda}_0(t)= P[X_c\leq (w_+-w)] &=& \frac{\lambda_0}{\lambda_0+n_X(t)a+m_Y(t)b+(n_{Xc}(t)+m_{Yc}(t))}\\
	&\stackrel{by\ (\ref{combined_Breslow})}{=}& \frac{1}{[n_{X}(t)+m_Y(t)]\hat{c}_{PL}+(n_{Xc}(t)+m_{Yc}(t))}.
\end{eqnarray*}
Solve the equation for $\lambda_0$. It follows
\[ \lambda_0=\frac{n_X(t)a+m_Y(t)b+(n_{Xc}(t)+m_{Yc}(t))}{[n_{X}(t)+m_Y(t)]\hat{c}_{PL}+(n_{Xc}(t)+m_{Yc}(t))-1}. \]
Consider that all patients in Example 4 have exponential survival time of various intensity ($a, b$ and 1 respectively), the law of large number guarantees
\[ \frac{n_X(t)}{n+m}\rightarrow pe^{-at}, \frac{m_Y(t)}{n+m}\rightarrow (1-p)e^{-bt}\ \mathrm{and}\ \frac{n_{Xc}(t)+m_{Yc}(t)}{n+m}\rightarrow e^{-t}. \]
Hence the limit of $\lambda_0$, the estimated hazard for a control group patient living within the infinitestimal time interval $(t,t+dt)$ is 
\[ \lambda_0\approx\frac{pe^{-at}a+(1-p)e^{-bt}b+e^{-t}}{(pe^{-at}+(1-p)e^{-bt})c^*_{PL}+e^{-t}}. \]
This is formula (\ref{combined_Breslow_limit}).

\subsection{Inequalities for different versions of the combined treatment effect}
We have investigated various definitions of the overall treatment effect for a mixed patient population. Here we are going to demonstrate the quantitative relationship between these definitions assuming the simplest modelling setup as was described in Example 1. That is, the treatment group indicator $Z\sim Bernoulli(q)$ is the only covariate for the proportional hazard models. The survival times of the enrolled patients follow proportional hazard models with hazard ratios $a=e^\alpha$ and $b=e^\beta$ respectively for trial 1 and 2. The ratio of the sample sizes of the two trials always equals $p:(1-p)$. We have the following candidates for the definition of the overall (log) hazard ratio: \\
i) $\theta_L=p\alpha+(1-p)\beta$,\\
ii) $c_L=pa+(1-p)b$,\\
iii) $c^*_{PL}$ as the solution to equation (\ref{Lin_solution_single_covariate_simple}),\\
iv) $c^*_{HM}$, the harmonic mean of $a$ and $b$ as was defined in  (\ref{c_hm}).\\

\noindent
\textit{\textbf{Proposition 3.} For arbitrary $a<b$, $q\in (0,1)$ and $p\in (0,1)$, the following inequalities always holds true:\\
1) $a<c^*_{HM}<\exp(\theta_L)<c_L<b$;\\
2) $a<c^*_{PL}<c_L<b$.\\
}

\noindent
\textit{Proof.} The harmonic mean $c^*_{HM}>a$ is trivial.\\

\noindent
To compare $c^*_{HM}$ with $\exp(\theta_L)$, consider their ratio
\[ \frac{\exp(\theta_L)}{c^*_{HM}}=p\left(\frac{b}{a}\right)^{1-p}+(1-p)\left(\frac{a}{b}\right)^p. \]
Denote the above ratio by $R(a,b)$ and calculate
\[ \frac{\partial R(a,b)}{\partial b}=p(1-p)\left(\frac{a}{b}\right)^p(1/a-1/b)>0, \forall a<b. \]
Combined with the fact that $R(a,b=a)=1$, it indicates $R(a,b)>1$ for any $b>a$ therefore we proved
\[ c^*_{HM}<\exp(\theta_L). \]
Using the Jensen's inequality with the convex function $g(x)=e^x$, it is easy to prove that
\[ \exp(\theta_L)=e^{pa+(1-p)b}<pe^a+(1-p)e^b=c_L. \]
It is also trivial to see that the algebraic mean $c_L<b$.\\

\noindent
For simplicity, denote the integrand in (\ref{Lin_solution_single_covariate_simple}) by $f(u,a,b,c)$. To compare $a$ and $c^*_{PL}$, note that
\[ f(u,a,b,c=a)=\frac{(1-q)e^{-u}+pqae^{-au}+(1-p)qbe^{-bu}}{(1-q)e^{-u}+pqae^{-au}+(1-p)qae^{-bu}}\cdot e^{-u}> e^{-u}, \forall a<b. \]
Hence $\int_0^\infty f(u,a,b,a)du>1$. To achieve an equality in (\ref{Lin_solution_single_covariate_simple}), one must have the solution $c^*_{PL}>a$ because $f(u,a,b,c)$ is decreasing with respect to $c$.\\

\noindent
To compare $c^*_{PL}$ with $c_L$, observe that
\begin{eqnarray*}
	&&[pqae^{-au}+(1-p)qbe^{-bu}]-[pqc_Le^{-au}+(1-p)qc_Le^{-bu}]\\
	&=& p(1-p)q(a-b)(e^{-au}-e^{-bu})<0
\end{eqnarray*}
Therefore $f(u,a,b,c_L)<e^{-u}$ and equivalently $\int_0^\infty f(u,a,b,c^*_{PL})du<1$. To achieve an equality in (\ref{Lin_solution_single_covariate_simple}), one must have $c^*_{PL}<c_L$. $\Box$\\

\noindent
One may also be tempted to find a fixed order for $c^*_{HM}$ vs $c^*_{PL}$ and $c^*_{PL}$ vs $\exp(\theta_L)$. However, close examination of the algebraic definitions implies that the result of these comparisons depend on the values of the parameters. The rule of thumb is, $c^*_{HM}<c^*_{PL}<\exp(\theta_L)$ for most $a<b<1$ of practical importance. Only when $a$ is extremely small (typically smaller than 0.2) one can observe $\exp(\theta_L)<c^*_{PL}$. When $1<a<b$, $c^*_{PL}$ and $c^*_{HM}$ are pretty close to each other and Proposition 3 indicates that $c^*_{HM}<\exp(\theta_L)$. Here is an example demonstrating the relationships of these estimators.\\

\noindent
\textbf{Example 6.} Let $p=0.5$, $q=0.5$. We plotted the percentage differences between the three estimators in Figure 2. The second estimator in the list is always the basis for comparison.\\

\begin{figure}[ht!]
\centering
\includegraphics[width=150mm]{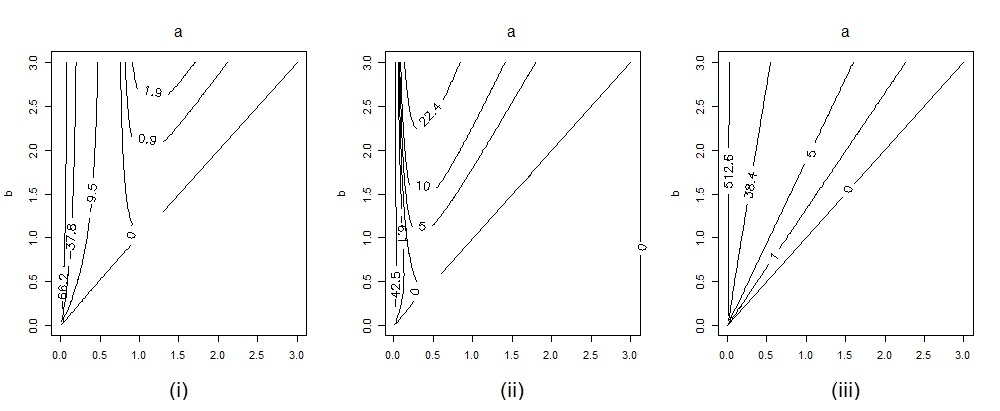}
\caption{Contour plots for the percentage differences between the estimators. (i) ${c}^*_{HM}$ vs ${c}^*_{PL}$. (ii) $\exp({\theta}^*_L)$ vs ${c}^*_{PL}$. (iii) $\exp({\theta}^*_L)$ vs ${c}^*_{HM}$.}
\end{figure}

\noindent
Note that when $a$ assumes a decent value, e.g., any number greater than 0.5, the linear log hazard ratio estimate $\hat{\theta}_L$ always leads to a conservative definition for the overall treatment effect compared the limit of the MPLE $c^*_{PL}$ or the limit of the MLE $c^*_{HM}$. The bias in $\exp({\theta}_L)$ can be higher than 38\% depending on the distance between $a$ and $b$. Table 1 are the exact values of the overall hazard ratios defined for various combinations of $a$, $b$:\\

\begin{table}[h]
\caption{Values of the overall hazard ratios}
\centering
\begin{tabular}{lll|llllll}
\hline
                       &                                           & b & 0.5 & 1.0 & 1.5 & 2.0 & 2.5 & 3.0 \\ \hline
\multicolumn{1}{l|}{}  & \multicolumn{1}{c|}{\multirow{4}{*}{0.5}} & $c^*_{HM}$ & $0.5\ \ $   & 0.662 & 0.741   & 0.792   & 0.823   & 0.847   \\
\multicolumn{1}{l|}{}  & \multicolumn{1}{c|}{}                     & $c^*_{PL}$ & $0.5\ \ $   & 0.682 & 0.781   & 0.848   & 0.892   & 0.925   \\
\multicolumn{1}{l|}{}  & \multicolumn{1}{c|}{}                     & $\exp(\theta_L)$ & $0.5\ \ $ & 0.705 & 0.857   & 0.994   & 1.107   & 1.216   \\
\multicolumn{1}{l|}{}  & \multicolumn{1}{c|}{}                     & $c^*_L$ & $0.5\ \ $   & 0.750 &   0.992 & 1.248   & 1.490   & 1.747   \\ \cline{2-9} 
\multicolumn{1}{l|}{}  & \multicolumn{1}{l|}{\multirow{4}{*}{1.0}} & $c^*_{HM}$ &    & 1.0 & 1.202   & 1.340   & 1.433   & 1.507   \\
\multicolumn{1}{l|}{} & \multicolumn{1}{l|}{}                     & $c^*_{PL}$ &    & 1.0 & 1.198   & 1.327   & 1.409   & 1.471   \\
\multicolumn{1}{l|}{a}  & \multicolumn{1}{l|}{}                     & $\exp(\theta_L)$ &    & 1.0 & 1.225   & 1.420   & 1.582   & 1.738   \\
\multicolumn{1}{l|}{}  & \multicolumn{1}{l|}{}                     & $c^*_L$ &    & 1.0 & 1.248   & 1.505   & 1.747   & 2.003   \\ \cline{2-9} 
\multicolumn{1}{l|}{}  & \multicolumn{1}{l|}{\multirow{4}{*}{2.0}} & $c^*_{HM}$ &    &  &    & 2.0   & 2.219   & 2.402   \\
\multicolumn{1}{l|}{}  & \multicolumn{1}{l|}{}                     & $c^*_{PL}$ &    &  &    & 2.0   & 2.212   & 2.375   \\
\multicolumn{1}{l|}{}  & \multicolumn{1}{l|}{}                     & $\exp(\theta_L)$ &    &  &    & 2.0   & 2.232   & 2.452   \\
\multicolumn{1}{l|}{}  & \multicolumn{1}{l|}{}                     & $c^*_L$ &    &  &    & 2.0   & 2.245   &   2.502
\end{tabular}
\end{table}

\def\bibindent{1em}

\end{document}